\documentstyle[12pt]{article}

\newtheorem{lemma}{Lemma}[section]

\newtheorem{proposition}[lemma]{Proposition}

\newtheorem{conjecture}[lemma]{Conjecture}
\newtheorem{definition}[lemma]{Definition}
\newtheorem{remark}[lemma]{Remark}

\newcommand{\selabel}[1]{\label{se:#1}}

\newcommand{\prlabel}[1]{\label{pr:#1}}
\newcommand{\prref}[1]{Proposition~\ref{pr:#1}}

\newcommand{\cjlabel}[1]{\label{cj:#1}}
\newcommand{\cjref}[1]{Conjecture~\ref{cj:#1}}
\newcommand{\relabel}[1]{\label{re:#1}}
\newcommand{\reref}[1]{Remark~\ref{re:#1}}

\newcommand{\delabel}[1]{\label{de:#1}}
\newcommand{\deref}[1]{Definition~\ref{de:#1}}

\newcommand{\eqref}[1]{(\ref{eq:#1})}

\newcommand{\Ker}{{\rm Ker}\,}

\newcommand{\End}{{\rm End}}

\newcommand{\sgn}{{\rm sgn}\,}
\newcommand{\Id}{{\rm Id}\,}

\def\ot{\otimes}

\def\cp{{\bf CP}^n}

\def\gggg{{\bf g}}
\def\hhhh{{\bf h}}
\def\nnnn{{\bf n}}

\def\real{{\bf R}}

\def\complex{{\bf C}}

\def\r#1{\mbox{(}\ref{#1}\mbox{)}}
\def\Uq{U_q(sl(n))}
\def\tri{\triangleright}
\def\uqsl{U_q(sl(2))}
\def\gl{\gamma_{\lambda}}

\def\ww{W^{\ot 2}}

\def\al{\alpha}
\def\la{\lambda}
\def\ml{m_{l}}
\def\lp{\wedge_+}
\def\lm{\wedge_-}
\def\lpm{\wedge_{\pm}}
\def\pp{P_+}
\def\pmm{P_-}
\def\ppm{P_{\pm}}
\def\ssl{s_{\lambda}}
\def\si{\sigma}

\def\De{\Delta}
\def\uq{U_q({\bf g})}

\def\h{{\hbar}}
\def\com{{\cal O}_{\omega}}

\def\hhhh{{\bf h}}

\def\nnnn{{\bf n}}

\def\bea{\begin{eqnarray}}
\def\eea{\end{eqnarray}}
\def\beq{\begin{equation}}          
\def\eeq{\end{equation}}

\def\aaaa{\cal A}
\def\aaa{\cal A}
\def\ah{{\cal A}_{\h}}
\def\adva{{\cal A}^{\ot 2}}
\def\ac{{\cal A}_{c}}
\def\ot{\otimes}
\def\om{\omega}

\def\Ml{M_{\lambda}}
\def\oMl{\overline M_{\lambda}}
\def\oMn{\overline M_{\nu}}
\def\End{{\rm End\, }}

\def\Ker{{\rm Ker\, }}
\def\Cas{{\cal C }}
\def\Cass{{\cal C }_{sl}}
\def\oCas{\overline{\cal C }}
\def\ev{{\rm ev}\,}
\def\root{{\rm root\, }}
\def\diag{{\rm diag\, }}
\def\Im{{\rm Im\, }}
\def\Id{{\rm Id }}

\def\id{{\rm Id\, }}
\def\Vect{{\rm Vect\,}}

\def\tr{{\rm tr\,}}
\def\det{{\rm det}}
\def\codet{{\rm codet}}
\def\dim{{\rm dim}\,}
\def\Gr{{\rm Gr}}
\def\udim{\underline{\rm dim}\,}
\def\dem{{\rm det}^{-1}}
\def\Fun{{\rm Fun\,}}
\def\ind{{\rm ind\,}}

\def\rank{{\rm rank\,}}
\def\bi{{\rm bi{\hbox{-}}rank\,}}
\def\vv{V^{\ot 2}}
\def\Tm{T^m(V)}
\def\Tn{T^n(V)}
\def\Tl{T^l(V)}
\def\W{{\bf V}}

\def\Wl{{\bf V}_{\la}}
\def\sgn{{\rm sgn}\,}
\def\res{{\rm res}\,}
\def\Sm{{\cal S}(m)}
\def\kSm{k[{\cal S}(m)]}
\def\S{{\cal S}}

\def\Vl{{\cal V}_{\lambda}}
\def\Vm{V_{\mu}}
\def\Vl{V_{\la}}
\def\Vn{V_{\nu}}
\def\clmn{c_{\la,\mu}^{\nu}}
\def\SM{$SL(n)$-Mod}
\def\rs{\rho_S}
\def\UV{{\cal SW}(V)}
\def\oVl{\overline{V}_{\la}}
\def\Tm{T^m(V)}
\def\de{\delta}

\def\R{{\cal R}}

\begin{document}
\title{Schur-Weyl Categories and non-quasiclassical Weyl type Formula}
\author{Dimitri Gurevich and Zakaria Mriss\\
Universit\'e de Valenciennes\\
Valenciennes, France}
\maketitle

\begin{abstract}\noindent
To  a vector space $V$ equipped with a non-quasiclassical involutary
solution of the
quantum Yang-Baxter equation and a partition $\la$,
we associate
a vector space $\Vl$ and compute its dimension. The functor $V\mapsto
\Vl$
is an
analogue of
the well-known Schur functor. The category generated by the objects
$\Vl$
is called the Schur-Weyl  category.
We suggest a way to construct some related twisted varieties looking
like
orbits of semisimple elements in $sl(n)^*$. We consider in detail
a particular case of
such ``twisted orbits", namely the twisted non-quasiclassical
hyperboloid and
we define
the twisted Casimir operator on it. In this case, we obtain a formula
looking like the  Weyl formula, and describing the asymptotic behavior
 of the function $N(\la)=\{\sharp\,\la_i\leq\la\}$,
where $\la_i$ are the  eigenvalues of this operator.
\end{abstract}

\section{Introduction}
It is well-known that the motivation to introduce and develop the
theory of quantum groups arose from the theory of integrable system.
Quantum groups (QG)
and their dual objects  (``quantum cogroups") supplied us an adequate
language to describe symmetries of some integrable models. However,
these objects have found another important application.
It turned out that they
provided us with a natural way to enlarge the framework of classical
 geometry. Namely, it was recognized that the usual flip $\si$
occurring in numerous constructions of  commutative and
non-commutative\footnote{Note that
the term  ``non-commutative geometry" is now used abusively  in
different
senses. We
prefer  to reserve it for non-commutative geometry in the sense of
A.Connes.
This type geometry deals  with non-commutative algebras looking like
operator algebras
in non-twisted categories. These algebras are usually equipped with a
commutative trace
and an involution satisfying the classical property $(A\,B)^*=B^*\,
A^*$
and their
derivations are usually defined by means of the usual Leibniz rule.
Cyclic
(co)homology
is also defined by means of the classical flip $\sigma$.
In contrast, ``twisted" geometry
deals with algebras equipped with a twist different from the usual
flip.
Such algebras can be ``twisted commutative"
or ``twisted non-commutative". In the latter case they are realized as
operator algebras
in twisted categories in the spirit of \cite{DGK}.}
geometry can be replaced by some other twists, in particular, those
arising from the well known QG $\uq$.

By a twist we mean a solution of
the so-called quantum Yang-Baxter equation (QYBE)
$$S^{12}S^{23}S^{12}=S^{23}S^{12}S^{23}$$
where
\beq
S:\vv\to\vv \label{sym1}
\eeq
 is a linear operator, $V$ is a vector space over a basic field $k$
and as usual
$$S^{12}=S\ot\Id\qquad{\rm and}\qquad S^{23}=\Id\ot S.$$
  (In a more general context, $S^{12}$ and $S^{23}$
are treated as operators
acting on the tensor product of three vector spaces $U\ot V\ot W$.)

Thus, the twists arising from the QG $\uq$ are defined by
\beq
S=\si(\rho_U\ot \rho_U)\R, \label{dva}
\eeq
where $\R$ is the quantum universal R-matrix corresponding
to the QG in question and $\rho_U$ is a representation of $\uq$ into a
given  $\uq$-module $U$.

Let us mention two properties of this twist.
In the first place, it is a deformation of the flip $\si$.
The category whose twists possess this property as well as their
objects will be called {\em quasiclassical}.
Secondly, it is not involutory $(S^2\not=\Id)$ and it cannot be made
involutory
by a rescaling $S\to aS,\, a\in k$.

This second property gives rise to some difficulties, that do not
occur in
the case of an involutory twist (in the sequel called a {\em
symmetry}).
The following question  appears from the very beginning: which
algebras can be
considered as twisted analogues of commutative algebras and, in
particular,
which system
of equations compatible with the action of the QG $\uq$ gives rise to
a
``twisted
variety"? (We refer the reader to the paper
\cite{DGK} for a detailed discussion of this problem.)

Nevertheless, the fact that the twist \r{sym1} is quasiclassical gives
us a
criterion of a
``raison d'\^etre" for $\uq$-covariant algebras: such an algebra is of
interest if it is
a  flat deformation of its classical counterpart\footnote{Let us
recall that
a $k[[\h]]$-module $\ah$ is called
a flat deformation of a $k$-module $\aaa$ if $\ah$ and ${\aaa}[[\h]]$
are
isomorphic
as $k[[\h]]$-modules and
$\ah/\h\ah=\aaa$. Here $\h$ is a formal parameter, $k$ is the basic
field
 and ${\aaa}[[\h]]$ stands for the completion of ${\aaa}\ot k[[\h]]$
 in
the $\h$-adic
topology.}.

However, for  non-quasiclassical twists  this criterion is no
longer valid. Thus, it is not so evident what algebra arising from
such a
twist can be considered as
a ``twisted variety". In
the present paper we suggest a way to
construct some twisted non-quasiclassical varieties looking like the
orbits
of semisimple elements in $sl(n)^*$. (By abusing the language we use
the
term ``variety"
for
 the corresponding function algebra.)
Nevertheless,  we restrict ourselves to algebras connected with
symmetries
 (a way to generalize our scheme to some non-quasiclassical twists of
 Hecke type is discussed in the last
section).

To construct such a twisted variety we need first a tensor category
possessing a
``sufficiently large" supply of objects. We construct such categories
looking like that
of
$SL(n)$-modules and we call them
 Schur-Weyl (SW) categories. Hopefully, such a category can also be
 treated as
the one consisting of
$sl(V_S)$-modules, where $sl(V_S)$ is a twisted non-quasiclassical
analogue
of the Lie
algebra
$sl(n)$. However, these categories can be introduced directly without
any
(usual or
twisted) Hopf structure or twisted Lie algebra. The objects of such a
category
are labeled by  Young diagrams (up to some identification)
in a way similar to the classical one  but
their dimensions are different from the classical ones.

This implies some drastic modifications in the well-known asymptotic
Weyl
formula. In
its classical form it says that the function
\beq
N(\la)=\{\sharp\, \la_i\leq \la\}
\label{Nl}
\eeq
where $\la_n$ are eigenvalues of the Laplace-Beltrami operator on a
smooth
compact (pseudo)Riemannian
variety $M$ has the following asymptotic behavior
\beq
N(\la)\sim c \la^{n/2},\qquad n=\dim M \label{Wf}
\eeq
with some constant $c$ depending on the volume of $M$.

However,  this formula is no longer valid in a non-quasiclassical
case. In
the present paper
we  show that on a twisted non-quasiclassical hyperboloid (which is
the
simplest example of
a ``twisted non-quasiclassical orbit") the function \r{Nl} of the
twisted
Casimir
operator has an  exponential growth w.r.t. $\sqrt\la$. (Let us observe
that
on  symmetric
 orbits in $\gggg^*$ the Casimir operator is
equal up to a factor to the Laplace-Beltrami one if the latter is
$\gggg$-invariant.)

To explain the reason of this phenomenon let us recall first some
aspects
of ``twisted
linear algebra"  developed
essentially in \cite{G} and \cite{L}. Let us fix a symmetry \r{sym1}
and
associate to it symmetric
and skew-symmetric algebras in a natural way
\beq
{\rm Sym}(V)=\lp(V)=T(V)/\{\Im (\Id-S)\},\quad  \lm(V)=T(V)/\{\Im
(\Id+S)\}
\label{var1}
\eeq
where $T(V)$ is the free tensor algebra of the space $V$
and $\{I\}$ stands for the ideal generated by  a subset $I\subset
T(V)$.
(These algebras are also well defined for the so-called
{\em Hecke symmetries}, i.e.
twists satisfying the relation
$$(q\,\Id -S)(\Id+S)=0,\qquad q\in k$$
if $q$ is generic.)

Let us remark that the algebra $\lp(v)$ is S-commutative (or simply
commutative)
in the following sense.
We say that an algebra $\aaaa$ equipped with
a symmetry $S:\ \adva\to\adva$ is commutative if
\beq
\mu=\mu S \qquad {\rm and}\qquad  S\mu^{12}= \mu^{23}S^{12}S^{23}
\label{sym}
\eeq
where $\mu:\adva\to\aaaa$ is the product in $\aaaa$
(the second relation  implies a similar one with interchanged couples
of
indexes 12  and 23).

So, identifying as usual
\beq
\Fun(V^*)\approx\lp(V) \label{alg}
\eeq
we can treat $V^*$ as an  example of a twisted variety which is not,
however interesting
from the geometrical viewpoint. Nevertheless, we want to point out the
main
peculiarity of this
``variety": the supply of elements of the  algebra \r{alg} differs
drastically from
that in the classical case. The very useful tool allowing us to
measure
this supply
is the so-called  {\em Poincar\'e (or Poincar\'e-Hilbert) series}.
These series are defined for symmetric $\lp(V)$ and skew-symmetric
$\lm(V)$
algebras by
$$\ppm(t)=\sum \dim \lpm^l(V)t^l$$
where $\lpm^l(V)$ is the degree $l$ homogeneous component of the
algebra in question. However, the classical relation
\beq
P_+(t)P_-(-t)=\Id\label{pp}
\eeq
is valid for any symmetry (or a Hecke symmetry for a generic $q$, cf.
\cite{G}).

As shown in \cite{G}, there exist a lot of
symmetries \r{sym1} with  $\dim V=n\geq 3$ such that
the corresponding Poincar\'e series $\pmm(t)$ is a monic polynomial of
degree
$p<n$. The degree $p$ is called {\em rank} of the space $V_S$
and denoted  $\rank V=\rank V_S$.
(We use the notation $V_S$ for a vector space $V$ equipped with a
twist
\r{sym1}.)

Let us remark that in the classical case ($S=\si$)
we have $\pmm(t)=(1+t)^n$ and consequently,  $\rank V_S=\dim V_S$.
This is also valid for any symmetry (or a Hecke symmetry) being a
deformation of
the classical
flip.
A (Hecke) symmetry  and corresponding vector space $V_S$ whose
Poincar\'e series $\pmm(t)$ is a monic polynomial will be
called {\em even}.

Thus, even symmetries whose Poincar\'e series $\pmm(t)$ are  different
from
$(1+t)^n,\,\, n=\dim V$ cannot be obtained by a deformation of the
flip
$\sigma$. We
call them {\em  non-quasiclassical}. Other varieties, different from
$V^*$,
and more interesting from a geometric point, are discussed in this
paper.
They are also non-quasiclassical, i.e., they are not a deformation of
a
classical
variety, and the ``supply of elements" in the corresponding
algebras differs drastically from the classical one.
As a measure of this supply we consider the function \r{Nl}
corresponding
to the
``twisted Casimir operator". (Let us emphasize that this operator
arises from a
Casimir element, which is not analogue of the Casimir element from
$\uq$,
see Section
4).

This paper is organized as follows.
In Section 1, we recall some facts from \cite{G},
about a possible form of an even symmetry $\r{sym1}$. The aim is to
show that
the family of non-quasiclassical symmetries whose so-called
determinant is
central is
big enough. In Section 2, we introduce the Schur-Weyl
category generated by a vector space $V_S$ equipped with such a
symmetry.
This category
 is formed by
objects $\Vl$ arising from twisted analogues of the Schur functor and
their
direct sums.
The main result of this Section is the computation of $\dim \Vl$.
In Section 3 we introduce a twisted Lie algebra of $sl(n)$ type which
plays
the role of
symmetries of the SW category and define a twisted analogue of the
Casimir
element.
In section 4 we
consider a particular example of a twisted non-quasiclassical variety,
namely the
twisted non-quasiclassical hyperboloid, and give an estimation of the
function $N(\la)$
corresponding to the Casimir operator on it. We
end the paper with a discussion of those aspects of our approach that
can be
generalized
to  Hecke symmetries, and we suggest a possible way to obtain other
``twisted non-quasiclassical orbits".

Throughout the whole of the paper the basic field $k$ is $\complex$ or
$\real$.

{\bf Acknowledgment} One of the authors (D.G.) was supported by the
grant 
CNRS PICS-608.

\section{Even symmetries}\selabel{2}
We begin with the following observation. Usually, we assume that a
tensor
category (or a quasitensor category, in the terminology of \cite{CP})
is given, for example the category of $\uq$-modules), and we study the
properties of the objects of such a category. Our approach will be
completely
opposite: We begin with a basic object
$V_S$ and generate some tensor category from it.

Moreover, this category can be constructed without
any Hopf or twisted Hopf algebra. We use a twisted Lie algebra (its
enveloping algebra can be
equipped with a twisted Hopf algebra structure) only to describe
twisted orbits and define twisted Casimir operator.
 From this viewpoint, Hopf or twisted Hopf algebras
are derived objects themselves and they can be found from the
reconstruction theorems
(cf. \cite{M}),
although their explicit description is not always easy
(cf. \cite{AG1} where an attempt is made to describe
an analogue of the QG $\uq$ for some non-quasiclassical Hecke
symmetries).

Let us pass now to describing a possible form of an even symmetry.
Let us fix a space  $V=V_S$.
Let $T^m(V)=V^{\ot m}$ be the m-th tensor power of the space $V$ (with
$T^0(V)=k$) and
$T(V)=\bigoplus_{m=0}^{\infty} T^m(V)$ its free tensor algebra. The
symmetry $S$ can be
naturally extended to the tensor algebra: we have
$$S:\ T^m(V)\ot T^n(V)\to T^n(V)\ot T^m(V)$$
and therefore  $S:\ T(V)^{\ot 2}\to T(V)^{\ot 2}$ (we keep the
notation $S$
for the
extended  symmetry). Moreover, we assume that
$$S(a\ot x)=x\ot a \qquad \forall\, x\in \Tm, \qquad \forall\, a\in
k.$$

In fact, a symmetry allows us to equip the space $T^m(V)$ with a
representation
of the symmetric group $\S(m)$ in the natural way. We associate to an
elementary transposition
 $s^{i\,i+1}\in \S(m)$
the operator
$$S^{i\, i+1}=\Id_{i-1}\ot S\ot \Id_{m-i-2}$$
where $\Id_i$ is the identity operator on $T^i(V)$. Any element of the
symmetric group can be
expressed as a monomial of the elementary transpositions. By
substituting
in this monomial
the operators  $S^{i\,i+1}$ we get  a representation of the symmetric
group
$\Sm$ into the space
$\Tm$. Consequently, we have a representation of the group algebra
$k[\S(m)]$. It will be denoted
$\rho_S$.
(Note that we treat the space $T^m(V)$ as a left $k[\S(m)]$-module,
 i.e., $(AB)x=A(Bx)$ for any $x\in\Tm$ and $A,B\in \kSm$.)

\begin{remark}\relabel{2.1}\rm
This representation of the symmetric group in tensor powers of a
linear
space  has a very particular property:
the operators $S^{i\,i+1}$ and $S^{i+1\,i+2}$ are related by the
formula
$$S^{i+1\,i+2}=\si^{i\,i+1}\si^{i+1\,i+2} S^{i\,i+1}
\si^{i+1\,i+2}\si^{i\,i+1}\qquad 1\leq i\leq m-2. $$
More precisely, the space $\Tm$ is equipped with an action of
$\Sm\times\Sm$, one copy
of $\Sm$ being represented by $\rho_S$ and the other one by
$\rho_{\si}$
and they are
related by the above formula.
\end{remark}

Let us consider the corresponding symmetric $\lp(V)$ and
skew-symmetric
$\lm(V)$
algebras defined in
the Introduction, and the corresponding Poincar\'e series $\ppm(t)$.

\begin{definition}\delabel{2.2}
We say that a symmetry \r{sym1} or the corresponding space
$V=V_S$ is even (resp., odd) if the Poincar\'e series
$\pmm(t)$ (resp., $\pp(t)$) is a monic polynomial
(i.e. a polynomial with leading coefficient 1). For an even symmetry
$S$,
we call the degree of the polynomial $\pmm(t)$ the rank of $V$, and we
denote this by $\rank V$.
\end{definition}

\begin{remark}\relabel{2.3}\rm
As shown in \cite{P}, the Poincar\'e series of a Hecke
symmetry is a rational function (a proof in the case of symmetries
also
appeared in \cite{D}).
Let us assume that $\pmm(t)$ is a rational function with monic
numerator and  denominator, and no common factors in numerator and
denominator.
We introduce the bi-rank $\bi V_S=(p,q)$ as the ordered pair
consisting of
the degrees of the numerator and denominator of $\pmm(t)$.
Thus, for any even (resp., odd) symmetry, we have $\bi V_S=(p,\, 0)$
(resp., $\bi
V_S=(0,q)$). Let us remark that the notion of  bi-rank is a
generalization of
super-dimension (see also Remark 1.5).
\end{remark}

In the sequel, we deal with even symmetries. Our next aim is to
introduce  the dual space $V^*$.
A space $V^*$ is called {\em right dual} if there exists
an extension of $S$ to
\beq
(V\oplus V^*)^{\ot 2}\to (V\oplus V^*)^{\ot 2}\label{prolon}
\eeq
and an invariant pairing
$$V\ot V^*\to k.$$
``Invariant" means that this pairing commutes with $S$ in the
following sense
$$<\,\,,\,\,>^{12}S^{23}S^{12}=S<\,\,,\,\,>^{23}$$
where
$$<\,\,,\,\,>^{12}=<\,\,,\,\,>\ot \Id,\qquad
<\,\,,\,\,>^{23}=\Id\ot <\,\,,\,\,>.$$
Both sides of this formula are treated as operators acting on $V\ot
V\ot V^*$.
Hereafter we index the operators in question from left to right
(for example, in the above formula, the operator
$S^{12}$ acts on $\vv$ and $S^{23}$ acts on $V\ot V^*$).

In a similar sense we will speak about invariance of other linear maps
$V^{\ot i}\to V^{\ot j}$.

Let us show that for any even (Hecke) symmetry the right dual space
exists.

To an even symmetry $S$, we associate the projector
$$P_-^p:\ T^p(V)\to \wedge_-^p(V).$$
In view of \deref{2.2},  we have $\dim\Im P_-^p=1$.

Fix a base
$$\{x_i\},\,1\leq i\leq n=\dim V,\,\,\,x_i\in V.$$
Then the  projector $P_-^p$ can be described as follows
\beq
P_-^p:\ x_{i_1}\ot x_{i_2}\ot...\ot x_{i_p}\to u_{i_1 i_2...i_p}v
\eeq
where
\beq
v=v^{j_1 j_2...j_p}
x_{j_1}\ot x_{j_2}\ot...\ot x_{j_p}\in T^p(V) \quad {\rm and}
\quad u_{i_1 i_2...i_p}v^{i_1 i_2...i_p}=1.
\eeq
Any index that appears as an upper and a lower index is assumed to be
a summation index; we omit the summation symbol $\Sigma$.

\begin{definition}\delabel{2.4}
The element $v$ (resp., $u=u_{i_1 i_2...i_p}x^{i_1}\ot
x^{i_2}\ot...\ot
x^{i_p}$) is  called
determinant  (resp., codeterminant) and denoted $\det$  (resp,.
$\codet$).
\end{definition}

Let us remark that the couple $(\det,\,\,\codet)$ is defined up to a
change
$$det\to a\,\det,\,\, \codet\to a^{-1}\codet,\,\,a\not=0.$$
In the sequel, we assume that such a couple is fixed.

We say that the space $W\subset \Tm$ is invariant if
$$S(W\ot V)\subset V\ot W.$$
It is not difficult to see that the subspace  $\Im P_-^p\subset
T^p(V)$ is
invariant.
This follows from the fact that the projector
$P_-^p$ can be expressed as a polynomial in $S^{12},\,S^{23},...,
S^{p-1\,p}$.
However, the determinant itself is in general not central.
Remark that the determinant and codeterminant become simultaneously
central.

We introduce two operators $M=(M_i^j)$ and $N=(N_i^j)$ acting on $V$,
defined
by their matrices with respect to the basis $\{x_i\}$:
$$M_i^j=u_{i_1 i_2...i_{p-1}i}v^{ji_1 i_2...i_{p-1}}\quad {\rm and}
\quad
N_i^j=u_{ii_1 i_2...i_{p-1}}v^{i_1 i_2...i_{p-1}j}$$
(we use the notation of \cite{G}).
In \cite{G} it is shown that
\beq
M_i^j N_j^k= p^{-2}\de_i^k \label{prod}
\eeq
($\de_i^k$ is the Kronecker symbol), and
\beq
S(v\ot x_i)=(-1)^{p-1}pM_i^j (x_j\ot v),
\quad S(x_i\ot v)=(-1)^{p-1}p N_i^j (v\ot
x_j).\label{com}
\eeq

Let $\dem$ be a new formal generator such that
$$\dem\,\,\det=\det\,\,\dem=1.$$
This implies that
the commutation rule of $\dem$ with the elements of the space $V$ is
the
inverse of
\r{com}, namely
\begin{eqnarray*}
S(\dem\ot x_i)&=&(-1)^{p-1}pN_i^j (x_j\ot \dem)\\
S(x_i\ot \dem)&=&(-1)^{p-1}p M_i^j (\dem\ot x_j).
\end{eqnarray*}
Let us define the dual space $V^*$ by fixing a base $\{x^j\},\,1\leq
j\leq n$
such that $x^j$ is identified with
\beq
\dem\ot (v^{j_1 j_2...j_{p-1}j}x_{j_1}\ot x_{j_2}\ot...\ot
x_{j_{p-1}}).\label{det}
\eeq

We leave it to the reader to verify that  the pairing
$$<\,\,,\,\,>:V\ot V^*\to k,\quad  <\,\,,\,\,> x_i\ot
x^j\mapsto\de_i^j$$
is invariant; it suffices to verify that the element
$x^i\ot x_i$ is invariant, i.e.
$$S(y\ot (x^i\ot x_i))= (x^i\ot x_i)\ot y\qquad \forall y\in V.$$
This completes construction of the right dual space.

\begin{remark}\relabel{2.5}\rm
There exists another base $\{y^j\}$ of the space $V^*$,
such that the pairing
$$<\,\,,\,\,>:V^*\ot V\to k,\quad  <\,\,,\,\,>y^j\ot x_i
\mapsto\de_i^j$$
is invariant. Thus, the space $V^*$ can also be treated as the
left dual space of $V$, if we equip it with an appropriate pairing.

Such a base can be introduced as follows: let $y^i=C^i_j x^i$, where
$C^i_j =T_{ik}^{jk}$ and the operator $T=(T_{ik}^{jl})$ is defined by
$$S_{ij}^{kl}T_{km}^{in}=\de_m^l\de_j^n.$$
We say that a twist $S$ is ``invertible by column" if such an operator
$T$
exists
(it is
easy to see that this is independent of the choice of basis).
For any twist invertible by column the right and left dual spaces can
be
introduced and they can be identified. In particular, this  means that
the
extension \r{prolon} exists.

In the sequel, we also need the operator defined by the matrix
$B_i^j=T_{ki}^{kj}$, which is the inverse of the one defined by
$C_i^j$,
cf. \cite{G}.

Moreover, for any Hecke symmetry ($q$ is assumed to be generic),
there exists a complex consisting of the terms $V^i\ot (V^*)^j$ and a
differential arising from the operator which is
inverse to the pairing between $V$ and $V^*$. Such a complex was
called
 in \cite{G} Koszul complex of the second kind.

If the cohomology of this complex
is one-dimensional, then its generator is called deteminant.
The above determinant is a
particular case  of the  latter one. However, up to now it is not
clear
whether any
(Hecke) symmetry invertible by column has a determinant. Hopefully,
the
Poincar\'e
series of a Hecke symmetry invertible by column is a rational function
with
monic
numerator and denominator (see \reref{2.3}).

In \cite{G}, this complex was used to show that, for any even Hecke
symmetry
with a generic $q$, the polynomial $\pmm(t)$ is reciprocal. The case
of
symmetries was
considered previously in \cite{L}.
\end{remark}

\begin{definition}\delabel{2.6}
We say that the determinant $\det$ is central if
$$
(-1)^{p-1}pM=\Id \quad {\rm and} \quad (-1)^{p-1}pN=\Id.$$
\end{definition}

In fact, in view of \r{prod}, the first relation implies the second
one and
vice
versa.

If det is central, then the dual space $V^*$ can be identified with
$\wedge^{p-1}(V)$
since in the formula \r{det} we can omit the factor $\dem$. In other
words,
if det is
central, then the map
\beq
\wedge^{p}(V)\to k,\qquad \det \mapsto 1\label{pair}
\eeq
is invariant.

Now we want to describe a family of even symmetries with
non-quasi\-clas\-si\-cal
Poincar\'e
polynomial $\pmm(t)$ and central determinant.

Let us begin with the simplest case $p=\rank V=2$.
In this case the polynomial $\pmm(t)$ is equal to $1+nt+t^2$ where
$n=\dim
V\geq 2$
(the case $n=2$ corresponds
to the quasiclassical case).
Then $S$ can be represented as
\beq
S_{ij}^{kl}=\de_i^k\de_j^l-2u_{ij}v^{kl} \label{form}
\eeq
with
\beq
u_{ij} v^{ij}=1. \label{odinn}
\eeq

Thus, the determinant in this case is $v=v^{kl}x_k x_l$.
If no confusion is possible, we omit the sign $\otimes$.
It is not difficult to see that if  $S$ is of the form \r{form}, then
the QYBE
for it is equivalent to the relation
$$u\,v\,u^t\,v^t={1\over 4}\,\Id$$
or, in a more detailed form,
\beq
u_{ij}\,v^{jk}\,u_{lk}\, v^{ml}={1\over 4}\,\de_i^m\label{odin}
\eeq
where $u^t$ is the transposed of $u$.

In  \cite{G}, a classification of all solutions $(u,v)$ of the
equations
\r{odinn}-\r{odin} is given,
including a more general case of Hecke symmetries\footnote{Some
symmetries of this type were discovered independently in \cite{DL}.}.
However,
without any classification, it is easy to see that the system
\r{odin}-\r{odinn}
possesses a large set of  skew-diagonal solutions. These are solutions
$(u,v)$ for which the only non-trivial elements appear on the
skew-diagonals
of $u$ and $v$, or
$$u_{ij}=0=v^{ij} \qquad {\rm  if}\qquad  i+j\not=n+1.$$
Some of them satisfy a complementary condition of centrality of
 the determinant. In the case $p=2$, this condition takes the
 following form
\beq
2u_{ij}v^{jk}=- \de_i^k\qquad {\rm in \,\, particular,}
\qquad  2u_{i\, n+1-i}v^{n+1-i,\,i}=-1,\,\, i=1,...,n\label{dva1}
\eeq
if $u$ and $v$ are skew-diagonal (in the latter formula no summation
is
assumed).

We leave it to the reader to describe the family of solutions of the
system
\r{odinn}-\r{dva1}. We restrict attention to the case $n=3$.
In this case the family of the couples $(u,v)$ satisfying the system
\r{odinn}-\r{dva1}
is parameterized by two indeterminates: if we choose
$u_{13}=a$ and $u_{22}=b$, then we have
$$u_{31}=-a/x,\,\,v^{13}=x/(2a),\,\,v^{22}=-1/2b,\,\,v^{31}=-1/2a$$
where $x$ is a solution of the equation $x+x^{-1}=3$.

Our next aim is to describe a way to construct even symmetries of rank
greater then 2
with central determinant.

First, assume that two symmetries
\beq S_1:\vv_1\to\vv_1 \qquad {\rm  and}\qquad S_2:\vv_2\to\vv_2
\label{dve} \eeq
are given. We present a procedure that allows us to construct a
symmetry
$S$ acting
on the space $(V_1\oplus V_2)^{\ot 2}$.
In \cite{MM}, such a procedure is called gluing; for Hecke symmetries
it was suggested earlier in \cite{G}. Let us assume that $S$
transposes
$V_1$ and $V_2$
by means of the usual flip $\si$. Then $S$ is a symmetry and we have
\beq
\ppm(t, V)= \ppm(t, V_1)\, \ppm(t, V_2).
\eeq
Therefore, if $V_1$ and $V_2$ are even, then $V$ is also even, and
$\rank V = \rank V_1+\rank V_2$.
Moreover, the matrices $p(-1)^pM$ and $p(-1)^pN$ related to the
symmetry
$S$ and
measuring  non-centrality of the  determinant connected to $S$
are equal to the tensor product of the corresponding
matrices related to $S_1$ and $S_2$.
The first statement is shown in \cite[Proposition 4.4]{G} and the
second one is
obvious.

This procedure enables us to construct a big family of even
non-quasi\-clas\-si\-cal symmetries
of higher rank and with central determinant by starting  with
symmetries of
rank 2 possessing
this property.
For example, if we take two symmetries
$$S_1:\ \vv_1\to\vv_1\qquad {\rm and}\qquad S_2:\ \vv_2\to\vv_2.$$
with Poincar\'e polynomials
$\pmm(t,V_i)=1+n_it+t^2$, where  $n_i=\dim V_i\geq 2,\,\,i=1,2$, then
$$\pmm(t, V)= 1+nt+(n_1n_2+2)t^2+nt^3+t^4,\qquad n=n_1+n_2.$$
Thus,  for a fixed $n=\dim V\geq 4$ we can construct symmetries
\r{sym1}
such that $\rank V_S =4$ and the middle coefficient of $\pmm(t)$
is equal to $a(n-a)+2,\,\,a=2,3,...,n-2$
(for other coefficients of $\pmm(t)$ there is no choice).

These examples show that the Poincar\'e polynomial $\pmm(t)$ of an
even
symmetry
is not determined by the couple
$(\dim V,\,\,\rank V) $.

Another way to construct even higher rank symmetries arises from Schur
functors
 described in the next Section; hopefully, the determinant
 corresponding to
the space
$\Vl$  is central as well.

\begin{remark}\relabel{2.7}\rm
If $S$ is an even symmetry, then $-S$ is an odd symmetry,
i.e., the series $\pp(t)$ is a
monic polynomial. Using the above gluing procedure
we can produce mixed symmetries, whose Poincar\'e series $\pmm(t)$ are
rational
functions with monic numerators and denominators. However,
there exist symmetries possessing such a type Poincar\'e series
$\pmm(t)$
and such that
the corresponding space $V$ cannot be split into a direct sum of an
even
subspace and an odd
one. The simplest example of such a symmetry is the following one
(\cite{L}):
$V=span(x,y)$, with
$$S(x\ot x)=x\ot x+by\ot y,\,\,
S(x\ot y)=y\ot x,\,\,S(y\ot y)=-y\ot y,\,\,b\in k$$
This space can be split into an even and an odd subspaces  iff $b=0$,
and
in this case it becomes a super-symmetry.
\end{remark}

\section{Schur-Weyl categories}\selabel{3}
Let
$$\la=(\la_1,
\la_2,...,\la_k),\qquad\la_1\geq\la_2\geq...\geq\la_k,\qquad
\vert\la\vert=
\la_1+...+\la_k=m$$
be  a partition
of an integer $m$. The corresponding Young diagram will be also
denoted by
$\la$.

We consider the right regular representation of $\kSm$, i.e.
equip the algebra $\kSm$ with the natural structure of a right
$\kSm$-module.
Then this algebra can be presented as
$$\kSm=\oplus [\Ml]^{\oplus n_{\la}},\qquad n_{\la}=\dim [\Ml]$$
where $[\Ml]$ is the class of the pairwise isomorphic
irreducible $\kSm$-module corresponding to the partition $\la$
and $\la$ runs over all partitions of the integer $m$. We keep the
notation
$\Ml$ for a
representative of this class. Thus, $\Ml$ is understood to be a
$\kSm$-module equipped
with an embedding $\Ml\hookrightarrow \kSm$.

In particular, such an embedding arises from the following procedure.
Let us
convert the diagram $\la$ into a tableau by arranging
the integers $1,...,m$ by columns.
This means that we put in the first
column the numbers $1,2,..,\la_1'$, in the second one
$\la_1'+1,\la_1'+2,...,\la_2'$ and so on,
where
$$\la'=(\la_1', \la_2',...,\la_{k'}')\qquad
\vert\la'\vert=\la_1'+...+\la_{k'}'= m$$
 is the partition dual to that $\la$. We assign to this tableau the
 Young
symmetrizor
$$p_{\la}=c_{\la} r_{\la}$$
where $r_{\la}$ (resp., $c_{\la}$) is the symmetrizor by lines
(resp., the skew-symmetrizor by columns). Let us generate by the
elements
$p_{\la}\in\kSm$ the right
$\kSm$-module, i.e., consider the set
\beq
p_{\la}\, q,\qquad \forall \, q\in\kSm.
\label{odno} \eeq
 This $\kSm$-module is just
a representative of the family $[\Ml]$. This module and all related
objects
will be called canonical
(we will explain this choice of ``canonical"  embedding
$\Ml\hookrightarrow
\kSm$ later).

To a $\kSm$-module  $\Ml\hookrightarrow\kSm$, we associate the space
$$V(\Ml)=\Vl=\Im \rs(\Ml)$$
with $\rs:\ \kSm\to\End(\Tm)$ as above. Thus, the space $\Vl$ is 
equipped with an embedding $\Vl\hookrightarrow\Tm$ depending on 
the embedding $\Ml\hookrightarrow\kSm$. Let $[\Vl]$ be the class of
all 
such spaces $\Vl$ embedded in $\Tm$ in one or another way.

The image of the canonical tableau is just the set
$$\Im\rs(p_{\la}\, q)\,\qquad \forall q\in \kSm$$

Let $\oMl$ denote the two-sided module in $\kSm$ generated by all
$\Ml$
(in other words, $\oMl=[\Ml]^{\oplus n_{\la}}$).
Its image in $\Tm$ will be denoted $\oVl$.
Thus, we have
$$\oVl=\Im\rs(q_1\,p_{\la}q_2),\qquad\forall\, q_1,\,q_2\in \kSm$$
In contrast to the space
$\Vl$, which depends on the chosen embedding $\Ml\hookrightarrow
\kSm$,
$\oVl$ depends only on $\la$.

\begin{proposition}\prlabel{3.1}
We have
$$\dim \oVl=\dim \Vl\,\dim \Ml.$$
\end{proposition}

{\bf Proof}
The statement follows immediately from \cite[Lemma 6.22]{FH}.

Let us consider the classical case $(S=\si)$ in more detail.

Let $k=l(\la)$ be the length of the partition $\la$, i.e., number
of lines in the corresponding diagram. It is obvious that
the space $\Vl$ is trivial if $l(\la)>n=\dim V$. It is well known that
if we equip the initial space $V$ with an action of  the groups
$GL(n)$
the spaces $\Vl$ become irreducible $GL(n)$-modules as well as their
products with
$$(\det g)^p,\,\, p\in {\bf Z},\,\, g\in GL(n).$$
The family of irreducible $SL(n)$-modules (considered up to
isomorphisms)
coincides with $\{[\Vl],\,\,l(\la)\leq n\}$, up to
the following identification.
If two partitions differ by a shift, i.e., $\mu=\la+a$ (this means
that
$\mu_i=\la_i+a,\,\, 1\leq i\leq n$) then the corresponding irreducible
$SL(n)$-modules
are identified. This is motivated by the fact that a column consisting
of
$n$ entries
(n-column for short)
corresponds to
the representation defined by the determinant and
is trivial for the group $SL(n)$.  

So, in the category of $SL(n)$-modules
we can always reduce $\la$ with $\la_n\not=0$  to  $\mu$ with
$\mu_n=0$ by
means
of this
identification. In what follows this operation will be called the
{\em reduction procedure}.

Dimensions of the spaces $\Vl$ can be found from the well known
formula
\beq
\dim \Vl=\Wl(1,1,...,1)/\W(1,1,...,1). \label{weyl}
\eeq
Here $\W(z_1,z_2,...,z_n)$ is the ordinary Vandermonde determinant in
$n$
indeterminates and
 $\Wl(z_1,z_2,...,z_n)$ is  the generalized Vandermonde determinant
corresponding to the partition $\la$, and defined as follows
$$
\left\vert\matrix{z_1^{\la_1+n-1}&z_2^{\la_1+n-1}&\ldots&z_n^{\la_1+n-
1}\cr
z_1^{\la_2+n-2}&z_2^{\la_2+n-2}&\ldots&z_n^{\la_2+n-2}\cr
\vdots&\vdots&\ddots&\vdots\cr
z_1^{\la_n}&z_2^{\la_n}&\ldots&z_n^{\la_n}\cr}\right\vert
$$
Observe that
$$\W(z_1,z_2,...,z_n)=\Wl(z_1,z_2,...,z_n)\quad {\rm if}\quad
\la=(0,0,...,0).$$
The  quotient $\Wl(z_1,z_2,...,z_n)/\W(z_1,z_2,...,z_n)$ is called
the Schur function (polynomial)
in $n$ indeterminates corresponding to the partition
$\la$ and is usually denoted $\ssl=\ssl(z_1,z_2,...,z_n)$.
We do not consider Schur
functions in an infinite number of indeterminates, which can be
defined as
limits of the
above ones for $n\to\infty$. Thus, in virtue of \r{weyl} we have
\beq
\dim \Vl=\ssl(1,1,...,1)\label{dim}
\eeq
Note that if $l(\la)=l<n$ we put in the above formula
$\la_{l+1}=...=\la_n=0$.

Let us consider now the category \SM\ of all finite dimensional
$SL(n)$-modules.
The classes $[\Vl]$ with $l(\la)< n$
form a base of this category, i.e. any object of \SM\ is isomorphic
to
a direct sum of irreducible $SL(n)$-modules
$[\Vl]$. Thus, the tensor product of any two objects of this category
is
determined by that
of two basic objects. The latter product is given by the formula
\beq
[\Vl]\cdot [V_{\mu}]=\clmn [\Vn]. \label{fusion}
\eeq
The coefficients $\clmn$ occurring in this decomposition
can be found by means of  the Littewood-Richardson rule.
However, it is necessary to keep in the mind that if for some $\nu$
entering this sum we have $l(\nu)> n$ the component $[V_{\nu}]$
disappears and
if  $l(\nu)= n$ we
 reduce $\nu$ as  above.

The algebra consisting of finite sums with integer coefficients of
formal
objects
$[\Vl]$ equipped with the product \r{fusion}
is called the {\em fusion ring} of the group $G=SL(n)$ (or $SU(n)$).

Now we pass to the general case $S\not=\si$.
Although we do not have any  object of Hopf algebra type
(it will be introduced in the next Section), we can introduce a
category
looking like
that of
$SL(n)$-modules directly. Its objects are finite sums of the spaces
$\Vl$.
In order to
treat this category as a
twisted tensor category, we have to explain how to
decompose the  tensor product
$$\Vl\ot\Vm \hookrightarrow
\Tm\ot\Tn=T^{m+n}(V)\qquad\vert\la\vert=m,\,\vert\mu\vert=n$$
 into a direct sum of $\Vn$. This can be done as follows.

We apply the set of operators $\rs({\oMn})$
to the product $\Vl\ot\Vm$. This defines a projection of this product
onto the component $\overline{\Vn}$.

Observe that the twist
\beq
S:\ \Vl\ot\Vm\to\Vm\ot\Vl
\label{trans}
\eeq
is well defined.
By definition, it is the restriction of the twist
$S:\ \Tl\ot\Tm\to \Tm\ot \Tl$, where
$\Vl$ (resp. $\Vm$) is embedded into $\Tl$ (resp., $\Tm$).
It is left as an exercise to the reader to show
that the image of $\Vl\ot\Vm$ belongs to  $\Vm\ot\Vl$.

It is easy to see that the formula
\r{fusion} is valid with the same coefficients as in the classical
case,
up to some
modifications. The role of $n=\dim V$ is played by $p=\rank V$. More
precisely, a
component
$[\Vn]$ occurring in the formula \r{fusion} is replaced by 0 if
$l(\nu)>p=\rank(V)$ and  it is reduced
as above if $l(\nu)=p$,
assuming that the determinant defined in the previous section
is central. Finally, we recover just the same
fusion ring as in the classical case but with $n$ replaced by $p$.
This fact has been already mentioned in the mathematical literature
(cf.
\cite{B}),
 even in more general situation related to Hecke symmetries.

However, dimensions of the spaces $\Vl$  and the corresponding
Clebsch-Gordan coefficients (which are defined if we fix some bases in
the
components
$\Vl$) are drastically different from the classical ones.
We will now calculate these dimensions.

Let $\beta_1, \,\beta_2,...,\beta_p$ be
the roots of the Poincar\'e polynomial $\pmm(t)$ corresponding to an
even
space $V_S$
such that $\rank V_S =p$ and let $\al_i=-\beta_i$. Then
$$\pmm(t)=\prod (t+\al_i )=\prod (\al_i t+1)$$
(in the latter equality we use the fact that this polynomial is
reciprocal).

The following proposition is a generalization of
the formula \r{dim}.

\begin{proposition}\prlabel{3.2}
Assuming that $l(\la)\leq p$ we have
$$\dim \Vl=\ssl(\al_1,...,\al_p).$$
\end{proposition}

This results immediately from \prref{3.1} and the following.

\begin{proposition}\prlabel{3.3} (\cite{L})
The multiplicity of the irreducible
$\kSm$-module
$[\Ml]$ related  to the partition $\la$ in the $\kSm$-module $\Tm$ is
equal to
$\ssl(\al_1,...,\al_p)$.
\end{proposition}

{\bf Proof}
Let $\chi_k$ be the character of the $\Sm$-module $T^m(V)$ (as above,
the
algebra $k[\Sm]$ is represented by $\rs$), $\chi^{\la}$ the character
of
$[\Ml]$ and $\eta_m$ the character of the trivial representation of
the
group $\Sm$.
 Then the multiplicity of the irreducible $\kSm$-module
$[\Ml]$ is measured by the following quantity
\begin{eqnarray*}
\langle\chi^{\la},\,\chi_m\rangle
&=&\langle\det(\eta_{\la_i-i+j}),\,\chi_m\rangle_{\Sm}\\
&=&
\sum\sgn\pi\\
&&\langle\eta_{\la_1-1+\pi(1)}...\eta_{\la_n-n+\pi(n)},\,\chi_m\rangle
_{\Sm}\\
&=&
\sum\sgn\pi\\
&&~~\langle\ind^{\Sm}_{\S(\la_1-1+\pi(1))\times...\times\S(\la_m-1+\pi
(m))}
\eta_{\la_1-1+\pi(1)}...\eta_{\la_n-n+\pi(n)},\,\chi_m\rangle_{\Sm}\\
&=&
\sum\sgn\pi\langle\eta_{\la_1-1+\pi(1)}...\eta_{\la_n-n+\pi(n)},\\
&&~~
\res_{\Sm}^{\S(\la_1-1+\pi(1))\times...\times\S(\la_m-1+\pi(m))}
\chi_m\rangle_{\S(\la_1-1+\pi(1))\times...\times\S(\la_m-1+\pi(m))}\\
&=&
\sum\sgn\pi\langle\eta_{\la_1-1+\pi(1)},\,\chi_{\la_1-1+\pi(1)}
\rangle_{\S(\la_1-1+\pi(1))}...\\
&&~~\langle\eta_{\la_m-1+\pi(m)},\chi_{\la_m-1+\pi(m)}\rangle_{\S(\la_
m-1+\pi(m)
)}\\
&=&
\det(\langle\eta_{\la_i-i+j},\chi_{\la_i-i+j}\rangle_{\S(\la_i-i+j)}).
\end{eqnarray*}

Here the pairing in question is $\Sm$-invariant and  the Frobenius
reciprocity is used.

Since $\langle\eta_k,\chi_k\rangle_{\S(k)}$ is the multiplicity of the
trivial module in $T^k(V)$, it is
equal to $\dim \wedge^k_+(V)$. Thus, we have
$$\langle\chi^{\la},\chi_m\rangle=\det
(h_{\la_i-i+j}(\al_1,...,\al_p))=\ssl(\al_1,...,\al_p)$$
where $h_k(x_1,...,x_p)$ are complete symmetric polynomials.
Here we use the relations
$$\ssl(x_1,...,x_p)=\det((h_{\la_i-i+j}(x_1,...,x_p))$$
 and
\beq
\langle\eta_{\la_i-i+j},\chi_{\la_i-i+j}\rangle_{\S(\la_i-i+j)}=
\dim \wedge^k_+(V)=h_k(\al_1,...,\al_p).
\label{symf}
\eeq
This completes the proof.

Let us observe that if two partitions $\la$ and $\mu$ such that
$l(\la)\leq
p,\,\, l(\mu)\leq p$
differ by a shift we have $\dim \Vl=\dim V_{\mu}$.

The previous result makes very plausible the following.

\begin{conjecture}\cjlabel{3.4}
Let $\root\pmm(t)$ denote the set of the roots of the
polynomial $\pmm(t)$ and $-\root\pmm(t)$ be the set of the opposite
numbers. Then
$$-\root\pmm(t,\Vl)=W_{\la}(-\root\pmm(t,V))$$
where
$$W_{\la}(z),\,\,z=\{z_1,z_2,...z_p\},\,\,z_i\in {\bf C},\,\,
z_1...z_p=1$$
 is defined in the following way.
To $z$ we associate the diagonal matrix
$$\diag(z_1,z_2,...z_p).$$ 
Then the matrix corresponding to $z$ in the
$SL(p)$-module $\Vl$ is also diagonal. The set of its diagonal
elements
is denoted $W_{\la}(z)$.
\end{conjecture}

\begin{remark}\relabel{3.5}\rm
Besides dimensions of the objects in any twisted rigid category
(i.e. a category closed with respect to the functor $V\mapsto V^*$),
there are
also the so-called inner (or quantum) dimensions, given by
$$\udim V= \tr \Id_V$$
where $\tr:\ \End(V,V)\to k$ is the trace which is
well defined in any rigid twisted tensor category
(see \cite{CP} and the Section 4) and $\Id_V:\ V\to V$ is the identity
operator
considered  as an element of $\End(V,V)$\footnote{The space
$\End(U,V)$ is
identified
with
$V\ot U^*$ (here $U^*$ is the left dual of $U$) and is called the
space of
(left)
inner morphisms from $U$ to $V$.}.
>From results of \cite{G}, it follows that
$$\udim V=\rank V$$
 for any even symmetry.
Indeed, in view of \cite[Proposition 2.12]{G} we have
$$\tr\,\Id=T_{ij}^{ij}=p$$
where $T_{km}^{in}$ is the operator mentioned in \reref{2.5}.

If \cjref{3.4} holds, then this implies that for any even symmetry
$$\udim \Vl=\ssl(1,1,...,1)$$ where
the unity is taken $p$ times. If $p=n$, then this is just the
classical
formula. Thus, the
inner dimension
of an even object depends only on the rank, i.e., on the degree of the
polynomial $\pmm(t)$, while
the ordinary dimension depends on the roots of this polynomial,
i.e., on the whole polynomial.
\end{remark}

\begin{definition}\delabel{3.6}
 Let $V_S$ be a vector space equipped with a symmetry
\r{sym1} such that  $\rank V=p$ and the
determinant is central. We call Schur-Weyl (SW) category and denote
$\UV$
the twisted symmetric category whose objects
are the spaces $\Vl,\,\, l(\la)\leq p$ and their direct sums and whose
morphisms of the
objects $\Vl$ are of two types.

The first type morphism is by definition a linear map of the form
$$\rs(p):\Tm\to \Tm,\qquad p\in\kSm.$$
Such morphisms give rise to a change of embedding of a given object
$\Vl\hookrightarrow\Tm$.
The second type morphism arises from
the reduction procedure as follows.

Let $\Ml$ be the right $\kSm$-module canonically embedded into $\kSm$
and
$\Vl$ be the corresponding
subspace of $\Tm$. If the diagram $\la$ contains a p-column we apply
the
map \r{pair} and kill it.
 (This is just the motivation of the above ``canonical"
arrangement.) The inverse linear map which is well defined is also a
morphism by
definition. A morphism of two direct sums of objects $\Vl$ is by
definition
a map being
a morphism on each component.
As usual we say that a morphism is an isomorphism if its inverse
exists and
is a
morphism as well.
\end{definition}

Remark that the map  $V\longmapsto \Vl$ is a twisted analogue of the
well-known
 {\em Schur  functor} defined in the case $S=\si$ (cf. \cite{FH}).

Let us also observe that for two different embeddings
$\Vl\hookrightarrow
\Tm$, there
exists a morphism sending one of them to the other one.

\begin{remark}\relabel{3.7}\rm
It is worth saying that the definitions of SW category and
 Schur functor can be
naturally generalized to Hecke  symmetries
(cf. \cite{P}).  However, the corresponding twists \r{trans} in a
particular case
$\la=\mu$ are not Hecke symmetries anymore.
\end{remark}

\begin{remark}\relabel{3.8}\rm
Note that, in the classical case (when  the above category is
just that of $SL(n)$-modules) there exists another way to introduce a
decomposition
of the product
$\Vl\ot\Vm$ into a direct sum of irreducible modules, using the notion
of a
highest
weight element with respect to a triangular decomposition
of the algebra $sl(n)$. So, we can study the above decomposition
without
any embedding
irreducible $SL(n)$-modules into tensor powers of the basic space $V$.
Unfortunately, in the general case
such an approach is not yet elaborated, it is not even clear what a
triangular decomposition of the corresponding twisted algebra
(considered
in the next
Section) should be. So, the only way to introduce a category looking
like
that of
$SL(n)$-modules is the Weyl type scheme developed above.
\end{remark}

\section{The twisted Lie algebra $sl(V_S)$ and the twis\-ted Casimir
operator}
Consider a vector space $V=V_S$ equipped with a symmetry $S$ that is
invertible by column, and that has therefor a well-defined
(say left) dual space $V^*$. Identifying $\End V$ with $V\ot V^*$,
we can extend the symmetry $S$ to (see \reref{2.5})
$$S=S_{\End V}:(\End V)^{\ot 2}\to (\End V)^{\ot 2}.$$
Here we treat the elements of
 $\End(V)$ as left (inner) morphisms, and
in this setting the space $V$ becomes a left $\End(V)$-module.

Moreover, $\End(V)$ can be equipped in a natural way by
a twisted (generalized or S-) Lie bracket as follows
$$[\,\,,\,\,]=\circ(\Id-S),\,\, S=S_{\End V} $$
where $\circ$ is the operator product in $\End(V)$.
The space $V$ equipped with such a bracket will be denoted by
$gl(V_S)$.

Now consider the S-trace (or simply the trace), defined on
the algebra $\End(V)$ by
$$\tr:\ \End V\to k,\qquad \tr=<\,\,,\,\,>S,\qquad\End V=V\ot V^*.$$
Remark that the trace is invariant and symmetric. Thus, we have
\beq\tr\,[\,\,,\,\,]=0.\label{tra}
\eeq

Let $e_i^j$ be the element of $\End V$ for which
$$e_i^j (x_k)=\de_k^j\, x_i,$$
i.e.,  we identify $e_i^j$ with $ x_i\ot x^j$. Then
$\tr(e_i^j)=C_i^j$, where the operator $C$ is defined as in
\reref{2.5}.
It follows from \r{tra} that the traceless elements of $\End(V)$ form
a subalgebra with respect to the twisted Lie bracket mentioned above.
This subalgebra will be denoted by $sl(V_S)$.

It is worthwhile to mention
the bracket $[\,\,,\,\,]$ can be expressed in terms
 of $C_i^j$ in the case where ${\rm rank}(V)=2$ (cf. \cite{G}).

The algebras $gl(V_S)$ and $sl(V_S)$ are particular cases of twisted
(generalized
or S-) Lie algebras defined as follows.

\begin{definition}\delabel{4.1} (\cite{G})
$\gggg=(V_S, [\,\,,\,\,]:\ \vv_S\to V_S)$ is called a twisted
(generalized or S-)Lie algebra
if the bracket $[\,\,,\,\,]$ is invariant and skew-symmetric and if
the following twisted analogue of the
Jacobi relation holds:
$$[\,\,,\,\,][\,\,,\,\,]^{12}(\Id+S^{12}S^{23}+S^{23}S^{12})=0.$$
\end{definition}

Assume that $V_S$ is equipped with an invariant and symmetric (resp.,
skew-symmetric)
pairing. Then we can introduce twisted Lie algebras of $so$ (resp.,
$sp$)
type  as
the subalgebra of $sl(V_S)$ consisting of elements preserving this
pairing.

The enveloping algebra of a twisted Lie algebra $\gggg$ is defined
in the following natural way
$$U(\gggg)=T(\gggg)/\{x_i\ot x_j-S(x_i\ot x_j)-[x_i, x_j]\}.$$
This enveloping algebra can be made into a cocommutative Hopf algebra
(see \cite{G}). The
comultiplication is given by the formula
$$\De x_i=x_i\ot 1+1\ot x_i.$$
There is a version of the PBW Theorem for the enveloping algebra
$U(\gggg)$.

\begin{proposition}\prlabel{4.2}
There exists a natural isomorphism
$$\wedge_+(\gggg)\cong \Gr\,U(\gggg)$$
where $\Gr\, U(\gggg)$ is the graded quadratic algebra associated to
the filtered algebra $U(\gggg)$.
\end{proposition}

{\bf Proof}
The algebra $\wedge_+(\gggg)$ is Kozsul (cf. \cite{BG} for the
definition).
It follows from the exactness of the Koszul
complex of the first kind from \cite{G}. Then by \cite{BG} we have the
result.

We say that a linear space $W$ is a
$\gggg$-module if there exists a twist $S:\ \ww\to\ww$ which can be
extended
to
$$S=S_{\End}:\ (\End W)^{\ot 2}\to (\End W)^{\ot 2}$$
and an invariant linear map $\rho:\ \gggg\to\End(W)$ such that
the operators associated to the elements $x_i$ via $\rho$ satisfy the
same
relations as
the elements $x_i$ themselves in the enveloping algebra.
The map $\rho$ is called a representation
of the algebra $\gggg$.
(It defines a representation of the algebra $U(\gggg)$ as well.)

Our next aim is to realize the category $\UV$ as that of
$\gggg$-modules (with
$\gggg=sl(V_S)$).
It is obvious that there is only one way to do this.
We have already defined the action of $\gggg$ on the base space
$V=V_S$.
We can extend this action to any tensor power of $V$ by means of the
above comultiplication. Observe that this action commutes with the
symmetry,
because the comultiplication is cocommutative. Therefore
this action commutes with any morphism
$$\rs(p),\,\, \forall\, p\in\kSm\quad \forall\, m.$$
This implies that all elements of $\gggg$ map any space
$\Vl\hookrightarrow \Tm$ into itself. Passing from one embedding to
another one
corresponds to passing from one representation of $\gggg$ to
an isomorphic representation.
Otherwise stated, we can say that the the first type morphisms commute
with
the action
of $\gggg$.

A similar statement holds for the morphisms of the second type. This
ensues
from the following result.

\begin{proposition}\prlabel{4.3}
The  twisted Lie algebra $\gggg=sl(V_S)$ maps the determinant $v$
into 0.
\end{proposition}

{\bf Proof} The statement results from the following Proposition.

\begin{proposition}\prlabel{4.4}
 The following formula holds
$$X(v)=p\, (\tr X)\, v,\quad X\in \gggg=sl(V_s).$$
\end{proposition}

{\bf Proof}
The map
$$X\ot x\to X(x),\,\, X\in\gggg,\,\,x\in V_S$$
is invariant and  $v\in \lm^p$, so we have
$$X(v)=X(v^{i_1 i_2\cdots i_p}x_{i_1}\ot x_{i_2}\ot...\ot x_{i_p})=Q\,
v^{i_1 i_2\cdots i_p}(X(x_{i_1}))\ot x_{i_2}\ot\cdots \ot x_{i_p}$$
with
$$Q=\Id-S^{12}+S^{23}S^{12}+\cdots+(-1)^{p-1}S^{p-1\,p}\cdots
S^{23}S^{12}.$$
It is easy to see that $X(v)\in\lm^p$.
Indeed,
$$v^{i_1 i_2\cdots i_p}(X(x_{i_1}))\ot x_{i_2}\ot\cdots\ot x_{i_p}\in
V\ot\lm^{p-1}(V)$$
and the operator $Q$ maps the space $V\ot\lm^{p-1}(V)$ into
$\lm^p(V)$.

Thus the element $X(v)$ does not change if
we apply the projection $P_-^p$ to it. Setting $X=a_i^je_j^i$,
we have
$$P_-^p(X(v))=p\, v^{i_1 i_2\cdots i_p}P_-^p(a_{i_1}^jx_j\ot
x_{i_2}\ot\cdots\ot x_{x_p})=$$
$$p\, v^{i_1 i_2\cdots i_p}a_{i_1}^j\, u_{j i_2\cdots i_p} v=p\,
C^{i_1}_j
a_{i_1}^j v$$
since $C^{i}_j=v^{i i_2\cdots i_p}\, u_{j i_2\cdots i_p}$ (cf.
\cite{G}).
The proof is complete after we observe that
$$\tr(X)=C^{i_1}_j a_{i_1}^j$$
because of  the equality $\tr\, e_i^j=C_i^j$.

\begin{conjecture}\cjlabel{4.4}
 All the $\gggg$-modules $\Vl$ are irreducible and any irreducible
 finite
dimensional $\gggg$-module is isomorphic to one of them. Moreover, a
linear map
 between two
objects of the category $\UV$ is a morphism in the above sense if and
only
if it is
a $sl(V_S)$-morphism, that is,
it commutes with the action of  $sl(V_S)$.
Thus, hopefully, the category $\UV$ can be treated as the category of
$sl(V_S)$- (or $U(sl(V_S)$-) modules.
\end{conjecture}

\begin{remark}\relabel{4.5}\rm
Beside the above twisted Hopf algebra $U(\gggg)$, the
category in question can be treated as the category of modules over a
usual
Hopf algebra $H$. Its dual Hopf algebra (quantum cogroup) $H^*$ has
been
constructed in
\cite{G}. An explicit description of the algebra
$H$ (which is also well defined for Hecke symmetries)
is not so easy (cf. \cite{AG1}). Let us mention also the papers
\cite{B}
and \cite{P}
where the algebra $H$ is considered.

Comparing these two Hopf algebras (the usual one and twisted one)
we want to emphasize that the twisted Hopf algebra $U(\gggg)$
is more  suitable for our aims because namely in terms of
this algebra we can describe tangent space of a twisted variety and
introduce
the twisted Casimir operator playing the role of the Laplace-Beltrami
operator
in our approach  (see \reref{5.1}).
\end{remark}

Now we equip the algebra $gl(V_S)$ with the pairing arising from the
trace,
namely,
$$\langle e_i^j,\,e_k^l\rangle=\de^j_k\tr e_i^l=\de^i_k\, C_i^j.$$
A direct computation shows that
the element $\id=e^i_i$ is orthogonal to the algebra $sl(V_S)$.
Moreover, the operator
$$gl(V_S)\to sl(V_S),\,\, e_i^j\mapsto
f^j_i=e_i^j-p^{-1}\,C_i^j\,\id,\,\,\, p=\rank(V_S)=\tr\, \id$$
is a projection onto the algebra $sl(V_S)$. The elements
$\{f^j_i,\,1\leq
i,\, j\leq
n\}$ generate this algebra but they are not free ($f^i_i=0$).

Now we define the (quadratic) Casimir element in the algebra
$U(gl(V_S))$
(resp., $U(sl(V_S)$)) as follows
$$\Cas=\Cas_{gl}=B_i^j\,e^i_l\,e^l_j\qquad
({\rm resp.,}\quad \Cass=B_i^j\,f^i_l\,f^l_j).$$
The operator $B=(B_i^j)$ is defined in \reref{2.5}.
These two Casimir elements are related by the formula
\beq
\Cass=\Cas-\Id\ot \Id/p. \label{twoC}
\eeq

It is easy to see that these elements are invariant. Now we will show
that
their images in
$\End(\Vl)$ are scalar (this also follows from \cjref{4.4}, but we
will not
make use of
it) and  compute the corresponding eigenvalues.
The operators arising from the Casimir elements $\Cas$ and $\Cass$
will be
called
Casimir operators and they will be denoted by the same letters.

Let us begin with the Casimir operator $\Cas$. We have
$$\Cas\tri e_k=B_i^j\, e^i_n\,e^n_j\tri e_k= B_i^j\, e^i_k\tri
e_j=B_i^i\,
e_k=p\,e_k.$$
The symbol $\tri$ stands  for the action of the operator in question
on an
element.

Applying this operator to the product $e_k\, e_l\in V_S^{\ot 2}$, we
obtain
$$\Cas\tri( e_k\, e_l)=(\Cas\tri e_k)\, e_l+e_k\,
(\Cas\tri e_l)+2\oCas\tri ( e_k\, e_l)$$
where $\oCas$, the so-called {\em split Casimir operator} is defined
by
the formula
$$\oCas \tri (e_k\, e_l)= \ev^{12} \,\ev^{34} \,S^{23}(B_i^j\,
e^i_n\,e^n_j\,e_k\,
e_l).$$
From here on, $\ev$ is the evaluation operator defined by $\ev (A\ot
x)=Ax$
where
$A$ is an  operator and $x$ is an element.

Using the properties of the tensor $B_i^j$ (cf. \cite[Section 1]{G}),
we can easily prove the following result.

\begin{proposition}\prlabel{4.6}
We have $\oCas\tri(e_k\, e_l)=S(e_k\, e_l)$ and therefore the formula
$$\Cas\vert_{V_S^{\ot 2}} =2\, p\, \Id+2S$$
holds.
\end{proposition}

This implies the following relation
$$\Cas\vert_{V_S^{\ot m}} =m\, p\, \Id+2\sum_{1\leq i< j\leq m}
S^{ij}.$$
Restricting the operator $\Cas$ to the component
$\Vl\hookrightarrow\Tm$,
we find
\beq
\Cas\vert_{\Vl}=m\,p\,\id+2\,Q_{\la},\qquad {\rm where} \qquad Q_{\la}
=Q\vert_{\Ml}\qquad {\rm and} \qquad
Q=\sum_{1\leq i< j\leq m} S^{ij}.\label{cas}
\eeq
Observe that the element $Q$ is central in the algebra $\kSm$. So,
being
applied
to $\Ml$ (we consider it as an operator $p\to Q\, p$) it becomes a
scalar
operator.
Thus, we have $Q_{\la}=\gl\,\Id$.
Since all summands of $Q$ are of the same cyclic type
$(m-2,1,0,...,0)$
 we have
\beq
\gl=(m^2-m)\chi_{\lambda}(C_{(m-2,1,0\cdots0)})/(2\,\dim
\Ml)\label{gam}
\eeq
where $C_{(i_1\,i_2\cdots i_l)}$ is the conjugacy class with the
cyclic
type $(i_1\,i_2\cdots i_l)$
and $\chi_{\lambda}$ is the character of $\Ml$.

Now we can conclude the following.

\begin{proposition}\prlabel{4.7}
The Casimir element $\Cas$ being applied to $\Vl$ as above becomes a
scalar  operator and
it is given by formula \r{cas} with $Q_{\lambda}=\gl\,\Id$ and
$\gl$ defined by \r{gam}.
\end{proposition}

Let us now consider two particular cases of the formula \r{cas}. If
$\la=(m,0,...,0)$ then
\beq
\Cas\vert_{\Vl}=m\, p\,\id+(m^2-m)\id=(m^2+m(p-1))\Id \label{casym}
\eeq
and if $\la=(1,1,...,1)$ ($m$ times) then
$$\Cas\vert_{\Vl}=m\, p\,\id-(m^2-m)\id. $$

By \r{twoC} we have
\beq
\Cass\vert_{\Vl}=\Cas\vert_{\Vl}-m^2/p\,\Id=(mp+2\gl-m^2/p)\,\Id.
\label{eigsl}
\eeq

We can conclude that the eigenvalues of the Casimir operators depend
$\rank V_S$, but not on $\dim V_S$.

\section{Non-quasiclassical hyperboloid and Weyl type formula}

The principle aim of this Section is to find a twisted
non-quasiclassical
analogue of the  asymptotic Weyl formula \r{Wf}.
The role of the Laplace-Beltrami operator will be played by the
Casimir
operator $\Cass$.  First, we will describe
a ``twisted non-quasiclassical variety", namely
the ``twisted non-quasiclassical hyperboloid".
The drastic difference between the behaviour of the function $N(\la)$
in the
classical case
and  the non-quasiclassical will be clear from this example. As usual,
we
assume
that the determinant is central.

First we observe that
 the space $sl(V_S)$ is itself an object $\Vl$ of the category $\UV$
corresponding
to the diagram
$\la=(2, 1^{p-2})$ where  $p=\rank V_S$.
If  $p=\rank V_S=2$, then the diagram corresponding to $sl(V_S)$ is
$\la=(2)$.
In this Section, we restrict attention to this case.

Let $\gggg=sl(V_S)$, and decompose the space $\gggg^{\ot 2}$
into a direct sum of objects $\Vl$ in the category $\UV$.
This sum contains three components $\Vl$ with
$$\la=(4),\qquad \la=(3,1),\qquad \la=(2,2).$$

If we carry out  the reduction procedure from above, we can reduce
the diagrams $\la=(3,1)$ and $\la=(2,2)$ to respectively
$\la=(2)$ and $\la=(0)$.
But instead of doing this, we consider
the symmetric algebra of $\gggg$ and impose some equations which are
compatible
with the  action of the twisted Lie algebra $sl(V_S)$.

Namely, we consider the following quotient algebra
$$\ac=T(\gggg)/\{f_i^j\,f_k^l-S(f_i^j\,f_k^l),\,\,B^j_if^i_k\,f^k_j
-c\},
\,\,c\in k.$$
If $c\not=0$, then the algebra $\ac$ is called
a (twisted non-quasiclassical) hyperboloid.
The algebra
${\aaa}_0$ is called a (twisted non-quasiclassical) cone.

It is not difficult to see that the latter equation
is compatible with the $sl(V_S)$ action. Thus, this twisted variety is
introduced
by means of
a unique equation in the symmetric algebra $\lp(sl(V_S))$ similar
to a classical hyperboloid (or cone).

Moreover, similarly to the classical case, it it possible to show that
the algebra
$\ac$ is
 a direct sum of the components
$$V_{(2 m)}\subset sl(V_S)^{\ot m} \qquad m=0,\,1,\,2,...$$
By using the results of the previous Section we are able
to estimate the function $N(\la)$ for the Casimir operator $\Cass$.
But first we want to realize this operator as a second order twisted
differential operator.
In order to do this we
will say some words on twisted differential operators on the
hyperboloid
in question. Connected to this is the paper \cite{GRR} where some
aspects of differential calculus arising from symmetries are
considered.

Recall that a twisted vector field (or S-vector field)
on a twisted commutative algebra $\aaa$ is
an operator $X:\ \aaa\to \aaa$ satisfying the
Leibniz rule, which is well known for an involutory twist $S$:
\beq
X(a\circ b)=X(a)\circ b+\circ\,\,\ev (X\ot\id)\,S\,(a\ot b),\,\, a, \,
b \in
\aaa.\label{Leib}
\eeq
If ${\aaa}=\lp(V)$, then we can identify the space $\Vect({\aaa})$ of
all
(left) vector fields with
${\aaa}\ot V^*$;
here $V^*$ is the (left) dual space, with action extended to the full
algebra
${\aaa}$ by means of the  Leibniz rule \r{Leib}).

We say that a vector field $X\in \Vect({\aaa})$ is a vector field on a
factor algebra
${\aaa}/\{I\}$ if
$$X(a)\in \{I\}\qquad \forall\, a\in \{I\}.$$
For $a\in {\aaa}$, we consider the operator $a\tri b=ab$. Operators of
the
form $X+a$, with $X\in \Vect({\aaa})$ and $a\in {\aaa}$ are called
first order differential operators on $\aaa$. In a similar way, we can
define differential operators of order $n$ on $\aaa$ and its
quotients.

Replacing $V$ by $\gggg$ in the previous example, we obtain a
definition of
the (left)
vector fields on the algebra $\lp(\gggg)$ and its quotients. A
particular
case of such
a vector field is given by those arising from the (left) adjoint
action of
the algebra
$sl(V_S)$ onto itself
$$X\mapsto {\rm ad}_X\qquad {\rm with}\qquad {\rm ad}_X\,Y=[X,\, Y].$$

Thus, the twisted Lie algebra $sl(V_S)$ is represented in
the algebra $\lp(\gggg)$ in two ways. The first one is given  via the
map
introduced in the previous Section and the second one is realized by
the
adjoint
action. It is worth saying that these two actions of the algebra
$sl(V_S)$ on
$\lp(\gggg)$ coincide. We do not need this statement in the sequel.
Let us
only observe
that by realizing the elements of $sl(V_S)$ as twisted vector fields
we can
treat the
operator $\Cass$ as a second order differential operator on the
algebra $\ac$.
Moreover, if $c\not=0$, then it is the unique second order operator
which
is $sl(V_S)$-invariant. In the case $c=0$ there exists another second
order
invariant
 operator, namely the square of the invariant vector field defined by
$X(f)=m\,f$
where $f$ is a degree $m$  element of ${\aaa}_0$.

\begin{remark}\relabel{5.1}\rm
Let us describe briefly  a way to introduce the tangent 
space on the hyperboloid in question. Denote $F_i^j$ the vector field
corresponding to
the element
$f_i^j\in sl(V_S)$.
It is possible to show that the vector fields $F_i^j$ generate the
space
$\Vect(\ac)$ as
a (left) $\ac$-module if $c\not=0$. Moreover, these vector fields
satisfy
the relation
$B_i^jf_k^iF_j^k=0$. So, it is natural to define
the tangent space on the twisted non-quasiclassical hyperboloid
as the quotient of the free left $\ac$-module generated by the formal
generators
$F_i^j$  such that $F_i^i=0$ over
the left submodule generated by the element $B_i^jf_k^iF_j^k$.
Moreover,
this  module (called tangent)
 is projective and the corresponding projection is
an $sl(V_S)$-morphism.

Let us emphasize that the operators coming from the Hopf
algebra $H$ mentioned in \reref{4.5} are rather useless for describing
the
tangent
space.

A similar situation takes place for a quantum
$\uqsl$-covariant hyperboloid. Its tangent space can be described in
terms of
``braided vectors fields"  which are completely different from
those coming from $U_q(sl(2))$, in spite of a tradition
assigning the meaning of vector fields
to the images of the elements $X,\,Y,\, H \in U_q(sl(2))$
(see \cite{AG2}).
\end{remark}

Now we assume that $n=\dim V_S>2$
(the case $n=2$ corresponds to the classical
hyperboloid).

\begin{proposition}\prlabel{5.2}
On the hyperboloid in question the eigenvalules $\la_l$ of the
Casimir
 operator $\Cass$ and
their multiplicities $\ml$ are
$$\la_l=(2l)^2/2+2l, \,\,\,\ml= (\al^{2l+1}_2-
\al^{2l+1}_1)/(\al_2-\al_1),\,\,\,
l=0,1,2,...$$
where $\al_i,\,\, i=1,\,2$ are the roots of the equation $P_-(-t)
=1-nt+t^2=0$.
\end{proposition}

{\bf Proof}
The result follows immediately from the formulae \r{eigsl} and
\r{symf}.

Let $\al_2 > \al_1$. Then we have the following.

\begin{proposition}\prlabel{5.3}
The function \r{Nl} possesses
the upper and low limits given by
$$\overline{\lim} N(\la)=\beta\,\al_2^{\sqrt{2\la}+2},\qquad
\underline{\lim} N(\la)=\beta\,\al_2^{\sqrt{2\la}}$$
where $\beta$ is a positive constant.
\end{proposition}

Let us emphasize two principle difference between the behavior of the
function
$N(\la)$ in the classical and non-quasiclassical cases.
First, in the non-quasiclassical case the
function
$N(\la)$ has an exponential growth
with respect to $\sqrt\la$ and, secondly, it does not have a limit but
only
an upper and a lower limit.

\begin{remark}\relabel{5.4}\rm
We are not able to give any estimation for the constant
$\beta$ in the spirit of the classical Weyl formula since we do not
know
any twisted
analogue of the notion of volume.
\end{remark}

It is interesting to compare this result with the analysis of the
spectrum
of an
``exotic harmonic oscillator" arising from non-quasiclassical
symmetries
introduced in
\cite{GRZ}.

Let us discuss now the case $p=\rank V_S >2$. In this case it is not
so
easy to find
a system of equations which would define
a  twisted non-quasiclassical variety. We restrict attention to the
``twisted
orbits" looking like the projective space ${\bf CP}^n$
embedded as an orbit ${\cal O}$ in $su(n)^*$ (this means that the
decomposition of the corresponding
algebra into a sum $\oplus \Vl$
looks like that of $\Fun ({\bf CP}^n)$).  To describe such a
twisted  orbit we should impose some system of equations on the space
$sl(V_S)$.  It is not difficult to guess their general form
using $sl(V_S)$-covariance of the system but
the problem is to find some factors occurring in this system (cf.
\cite{DGK} where
this problem is discussed w.r.t. to $\Uq$-covariant ``orbits").
To find such a system we
can use a scheme close to that considered in \cite{DGK}. Let us
describe
its classical
version.

Let $\gggg$ be a simple Lie algebra and $\hhhh+\nnnn_++\nnnn_-$ its
triangular
decomposition. Fix an element $\om\in\hhhh^*$ and extend it to $\gggg$
by
setting $\om(\nnnn_{\pm})=0$. Let
$\com$ be the $G$-orbit of $\om$, where $G$ acts on $\gggg^*$ by
$Ad^*$
action.
Consider the algebra $\aaa=\Fun(\com)$ defined as the restriction of
the
algebra
$\lp(\gggg)$ to the orbit $\com$.
Its quantization can be realized in the following way. We associate to
the
element
$\om$ some infinite dimensional $\gggg$-module $M_{\om}$ called the
(generalized) Verma
module - its construction is described in \cite{DGK}.
Let
$$\rho_{\om}:\ T(\gggg)\to\End(M_{\om})$$
 be the corresponding representation of the tensor algebra $T(\gggg)$.
 Set
$\rho_{\h}=\h\,\rho_{\om/\h}$.
Then the algebra
 $$\ah=T(\gggg)[\h]/\{\Ker\,\rho_{\h}\}
=\Im \rho_{\h}T(\gggg)[\h] \subset \End(V_{\om})[[\h]]$$
is a flat deformation of the initial algebra in the sense of footnote
2
(cf. \cite{DGK} for detail).
Let us remark that the quantum object is realized as an operator
algebra.

The  passage from the algebra $\ah$ to $\aaa$ is usually called
``dequantization" , and
can be used to find system of equations describing the orbit in
question
(compare to \cite{DGK}).

Unfortunately, for a twisted Lie algebra $\gggg$  corresponding to a
non-quasi\-clas\-si\-cal symmetry, it is not clear what its
(generalized) Verma module is. However,
we can suggest
some discrete analogue of this method dealing with finite dimensional
$\gggg$-modules.

Fix an object $\Vl$ of the Schur-Weyl category (in particular, the
generating space
$V=V_S$  itself) and set
$$V_l=\wedge^l_+(\Vl),\,\,l=1,2,...\qquad {\rm (thus,}\quad
V_1=\Vl).$$
Let
$$\rho_l:\ T(\gggg)\to\End(V_l)$$
be the representation of the twisted Lie algebra $\gggg=sl(V_S)$,
which is the extension of
the representation $T(\gggg)\to\End(V_S)$.
Consider the representation $l^{-1}\rho_l:\ T(\gggg)\to\End(V_l)$
(the passage from $\rho_1$
to $l^{-1}\rho_l$ is an analogue of the above passage from
$\rho_{\om}$ to
$\h\,\rho_{\om/\h}$).

Put
$$I_l=\Ker l^{-1}\rho_l T(\gggg),\qquad {\aaa}_l=T(\gggg)/\{I_l\}.$$
Hopefully, the algebras ${\aaa}_l$ converge to a commutative algebra
$\aaa$
which
is
considered as  a ``twisted orbit". The  system of equations describing
this
``orbit"
can be found from this limit (cf. \cite{DGK}).

For the ${\bf CP}^n$-type twisted orbits mentioned above, the
corresponding
``function algebra" should be decomposed into  a direct sum of
components
$$\Vl \qquad {\rm with} \qquad \la=(0),\,(2, 1^{p-2}),\,
(4,2^{p-2}),\,
(6,3^{p-2}),...$$
By using the results of the previous Section,
it is not difficult to obtain an estimation
of the function $N(\la)$ for the Casimir operator $\Cass$ in this
case.

Let us indicate now which aspects of the above theory can be
generalized to
the Hecke
symmetries. As we have already said the construction of a category
$\UV$
and that of
Schur functor have natural analogues in the case when the symmetry
\r{sym1}
is of
Hecke type since the representation theory  of the Hecke algebra for a
generic $q$
looks like that of the symmetric group. If a Hecke symmetry $S_q$ is a
deformation  of
a symmetry $S$, the dimensions of the spaces $\Vl$ arising from
corresponding Schur
functors are stable during the deformation $S\to S_q$.

The problem is to find a reasonable way to define corresponding
twisted non-quasiclassical varieties, to
introduce their tangent spaces and to define an analogue of the
Casimir
operator.
We hope to treat this problem elsewhere.  We refer the reader to the
paper
\cite{AG2} where the problem is solved for a quantum hyperboloid
related to
$U_q(sl(2))$.

We conclude by saying that, in order to construct ``reasonable"
twisted
varieties
 related to
non-involutary non-quasiclassical twists,
a criterion of flatness of deformation can be very useful.
By considering twists $S_q$ which are deformations
of a symmetry $S$, we should first define a twisted variety arising
from
the latter
symmetry (as we said in the Introduction, an involutory case is easier
to
study) and
then deform it to a variety related to the twist $S_q$. This scheme
looks
like the one
used in \cite{DGK} for introducing some $\Uq$-covariant algebras,
but a commutative algebra replaces as initial point is replaced
by an $S$-commutative one with an involutary $S$.



ISTV,
Universit\'e de Valenciennes,
F-59304 Valenciennes,
France\\
e-mail gurevich@univ-valenciennes.fr\\ and
zmriss@univ-valenciennes.fr
\end{document}